\documentclass[11pt]{article}
\usepackage{epsfig,amsmath,latexsym}
\usepackage{amsfonts}
\newtheorem{theorem}{Theorem}[section]

\def\slfrac#1#2{\hbox{\kern.1em %
 \raise.5ex\hbox{\the\scriptfont0 #1}\kern-.11em %
 /\kern-.15em\lower.25ex\hbox{\the\scriptfont0 #2}}}

\newcommand{\hsp}{\hspace*{\parindent}}

\newcommand{\eeq}{\end{equation}}
\newcommand{\beql}[1]{\begin{equation}\label{#1}}
\newcommand{\bsq}{{\vrule height .9ex width .8ex depth -.1ex }}

\newcommand{\NN}{{\Bbb N}}

\newcommand{\ZZ}{{\Bbb Z}}

\newcommand{\sP}{{\cal P}}

\newcommand{\sS}{{\cal S}}
\newcommand{\sT}{{\cal T}}
\newcommand{\sW}{{\cal W}}

\makeatletter
\def\@sect#1#2#3#4#5#6[#7]#8{\ifnum #2>\c@secnumdepth
     \def\@svsec{}\else
     \refstepcounter{#1}\edef\@svsec{\csname the#1\endcsname.\hskip .75em }\fi
     \@tempskipa #5\relax
      \ifdim \@tempskipa>\z@
        \begingroup #6\relax
          \@hangfrom{\hskip #3\relax\@svsec}{\interlinepenalty \@M #8\par}%
        \endgroup
       \csname #1mark\endcsname{#7}\addcontentsline
         {toc}{#1}{\ifnum #2>\c@secnumdepth \else
                      \protect\numberline{\csname the#1\endcsname}\fi
                    #7}\else
        \def\@svsechd{#6\hskip #3\@svsec #8\csname #1mark\endcsname
                      {#7}\addcontentsline
                           {toc}{#1}{\ifnum #2>\c@secnumdepth \else
                             \protect\numberline{\csname the#1\endcsname}\fi
                       #7}}\fi
     \@xsect{#5}}
\def\@begintheorem#1#2{\it \trivlist \item[\hskip \labelsep{\bf #1\ #2.}]}

\def\plain{plain}\ifx\fmtname\plain\csname fi\endcsname
     
     \input docstrip
     \preamble

     Do not distribute the stripped version of this file.
     The checksum in the header refers to the documented version.

     \endpreamble
     \generateFile{here.sty}{t}{\from{here.doc}{}}
     \endinput
\fi
\ifcat a\noexpand @\let\next\relax\else\def\next{%
    \documentstyle[here,doc]{article}\MakePercentIgnore}\fi\next
\ifx\@Hxfloat\@Hundef\else\expandafter\endinput\fi
\let\@Hxfloat\@xfloat
\def\@xfloat#1[{\@ifnextchar{H}{\@HHfloat{#1}[}{\@Hxfloat{#1}[}}
\def\@HHfloat#1[H]{%
\expandafter\let\csname end#1\endcsname\end@Hfloat
\vskip\intextsep\vbox\bgroup\def\@captype{#1}\parindent\z@
\ignorespaces}
\def\end@Hfloat{\egroup\vskip \intextsep}

\makeatother

\catcode`\@=11
\renewcommand{\section}{
        \setcounter{equation}{0}
        \@startsection {section}{1}{\z@}{-3.5ex plus -1ex minus
        -.2ex}{2.3ex plus .2ex}{\large\bf}%
        }
\catcode`@=12

\begin{document}
\begin{center}
{\Large {\bf Wild and Wooley Numbers}} 

\vspace{1.5\baselineskip}
{\em Jeffrey C. Lagarias} \\
\vspace*{.2\baselineskip}
University of Michigan  \\
Ann Arbor, MI  48109-1109 \\
{\tt lagarias@umich.edu} \\
\vspace*{2\baselineskip}
(May 1, 2005) \\
\vspace{1.5\baselineskip}
{\bf Abstract}
\end{center}
\noindent

The wild integer semigroup $\sW(\ZZ)$ consists
of the integers in the 
 multiplicative semigroup 
generated by $\{ \frac{3n+2}{2n+1}: n \ge 0 \}$
and $\frac{1}{2}$. 
The wild numbers are the irreducible
elements  in $\sW(\ZZ)$.
The paper presents evidence
that the wild numbers are  the set of  
all prime numbers excluding $3$. 
The smaller semigroup  $\sW_0(\ZZ)$ of
integers in the semigroup 
generated by  $\{ \frac{3n+2}{2n+1}: n \ge 0 \}$
are called the Wooley integer semigroup,
and its irreducible elements are called Wooley numbers.
This semigroup is shown to be recursive, and
various open problems are formulated about
Wooley integers.

\vspace{1.0\baselineskip}
\noindent
%
%
%
\setlength{\baselineskip}{1.0\baselineskip}

\section{INTRODUCTION}
\hsp
This paper studies  a problem in multiplicative
number theory originating from a weakened form
of the $3x+1$ problem. Let $\sW_0$ signify the multiplicative semigroup
generated by all rationals $\{ \frac{3n+2}{2n+1}: n \ge 0 \}$,
which is the set 
$\{ \frac{2}{1}, \frac{5}{3}, \frac{8}{5},\frac{11}{7}, ...\}.$
That is, $\sW_0$ is the set of all finite products of
the generators, allowing repetitions.
Let  $\sW$ signify the larger
multiplicative semigroup generated by
$\sW_0$ together with $\{ \frac{1}{2}\}$. 
We call $\sW$ the {\em wild semigroup}
and  $\sW_0$ the {\em Wooley semigroup},
respectively. 
The question we consider is:  Which integers belong to
these semigroups? 

The sets of integer elements
$\sW(\ZZ) := \sW \cap \ZZ$
and $ \sW_0(\ZZ):= \sW_0\cap \ZZ$
themselves form multiplicative semigroups, which we term the 
{\em wild integer semigroup} and {\em Wooley integer semigroup},
and we refer to  their members as  ``wild integers''
and ``Wooley integers,'' respectively.
We have  the immediate implication that each
Wooley integer is a wild integer, but the
converse need not hold. 
The Wooley semigroup $\sW_0$ is a semigroup without unit, 
whereas the
wild semigroup $\sW$ is a semigroup with unit, and
the semigroups 
$\sW_0(\ZZ)$ and $\sW(\ZZ)$ inherit these properties.
Our particular  choice of terminology is explained at the end of
the introduction.

An {\em irreducible element} of a commutative semigroup $\sT$
is one that cannot
be written as a product of two nonunits 
(noninvertible elements) in the semigroup 
(see Gilmer \cite[sec. 1.6]{Gi84}).
We call the irreducible elements of the wild integer semigroup
{\em wild numbers}; similarly we christen 
 the irreducible elements  of
the Wooley integer semigroup {\em Wooley numbers}.
Thus  the  wild numbers are a subset of the wild integers;
and the Wooley numbers are a subset of the Wooley
integers.

It is immediately evident that $2$ is both a wild number and
a Wooley number. It is also easy to show that $3$ is not a wild
number, hence not  a Wooley number. However the nature of 
other wild  numbers or Wooley numbers is less  evident.
The object of this paper is to determine  properties of
wild numbers and Wooley numbers.
It turns out that the Wooley numbers have a complicated
and not completely  understood structure; by comparison,  
the wild numbers have a reasonably simple description. 

The wild numbers and the Wooley numbers 
differ in some significant ways. An odd integer 
$w$ is in the wild integer semigroup
if and only if  there is a nonnegative integer  $j$ such
that $2^j w$ is in the Wooley integer semigroup. 
At the level
of irreducible elements, we infer that if $w$ is a wild
number, then for some $j$ $2^j w$ is a Wooley number;
if $2^j w$ is a Wooley number, however,  we cannot 
(currently) decide whether $w$ must be a wild number.
At first glance the  Wooley numbers are the simpler objects
from a computational perspective. In section 2 we show 
that there is
an effectively computable procedure for
deciding whether a given rational $r$ belongs to $\sW_0$.
This leads to  an effectively computable procedure to
determine whether an integer  is a Wooley integer
and, if so, whether it is a Wooley number. 

In contrast, it is not 
immediately apparent if there exists an
algorithm  for recognizing whether a
given integer is a wild integer. For
a general multiplicative semigroup generated by a 
recursive set of  rational numbers, 
it seems plausible  that the problem of recognizing which
integers belong to the semigroup is  sometimes an undecidable problem,
We refer the reader to  Dyson~\cite{Dy64} (and Taitslin \cite{Ta62})
for analogous undecidability results for commutative semigroups.
On the other hand, if one were able to 
characterize directly all members
of a semigroup, this could lead to  a decision procedure.
In the case at hand, we formulate a conjecture  that  would 
provide such a characterization, 
and thus lead to a  decision procedure.

There is  a strong connection between these problems
and a weakened form of
the $3x+1$ problem, which motivated their study.
This notorious problem (see \cite{La85} or \cite{Wi98}) 
is concerned with  the iteration of the
 function $T: \NN \to \NN$ defined by 
$$
T(x) = \left\{
\begin{array}{cl}
\frac{3x+1}{2}  & \mbox{if} ~~x \equiv 1~~ (\bmod ~ 2), \\
~~~ \\
\frac{x}{2} & \mbox{if}~~ x \equiv 0~~ (\bmod ~2 )~.
\end{array}
\right.
$$
The $3x+1$ conjecture asserts that
for each positive integer $n$ there exists an $N$ such that
$T^{(N)}(n)=1$, where $T^{(N)}= T \circ T \circ \cdots \circ T$ (N terms) 
is the $N$-fold iterate of $T$. It has been verified for all $n$ with
$n \le 10^{17}$,
but remains an open problem.

The weakened version of the $3x+1$ problem to which we alluded
earlier was
proposed by Herschel Farkas \cite{Fa04}. It reads as follows:

\paragraph{Weak $3X+1$ Conjecture.} 
{\em Consider the semigroup
$\sS := \sW^{-1} = \{ w^{-1}: w \in \sW\}$ that
is generated by $\{ \frac{2n+1}{3n+2}: n \ge 0\}$ together
with the integer $2$. Then $\sS$ contains every
positive integer.} \\

\noindent Farkas \cite{Fa04} observed that the truth of 
the $3x+1$ conjecture
implies the truth of the weak $3x+1$ conjecture, because
the steps of the $3x+1$ iteration process can be encoded as 
certain products of generators using the 
semigroup multiplication in $\sS$. However, there are 
products of generators in the semigroup $\sS$ that do not correspond to
the $3x+1$ iteration, so the Farkas conjecture
is potentially easier to resolve than the $3x+1$ problem. 

Investigation of  the weak $3x+1$
conjecture led to questions  about the wild integer semigroup 
$\sW(\ZZ)$ as a possible aid  in its proof. Conversely, the  weak 
$3x+1$ conjecture has very strong implications about wild
numbers, as explained in section 3,
that led to the formulation of  the following conjecture:

\paragraph{Wild Numbers Conjecture.}
{\em The wild numbers consist of  the set of  all prime numbers,
excluding  $3$.   Equivalently, the wild integer 
semigroup $\sW(\ZZ)$  consists of 
all positive integers $m$ not divisible by $3$.
} \\

\noindent
This paper studies properties of wild and Wooley integers
that bear on these questions.
In section 2  we study Wooley integers. 
We give an effectively computable algorithm for
recognizing Wooley integers and Wooley numbers.
Using this approach we show that $20$ is a Wooley number.
We also report on computations finding various Wooley integers.
We raise the question whether 
the Wooley integer semigroup $\sW_0(\ZZ)$ is a free commutative
semigroup  and present
evidence suggesting that it is not.
In contrast, in section 3 we  show that the weak
$3x+1$ conjecture implies that 
the wild integer semigroup $\sW(\ZZ)$ is a free 
commutative semigroup with unit.

In section 3  we begin by demonstrating  that there are infinitely many
wild numbers. Then we 
show that the weak $3x+1$ conjecture implies
strong restrictions on wild numbers---indeed, we prove
that it implies
the wild numbers conjecture. This lends strong
support to the conviction that the wild numbers
conjecture must be true, since the weak $3x+1$ 
conjecture itself would follow from the $3x+1$ conjecture,
for which there is extensive evidence.
We also deduce a
 converse assertion to the effect that 
the wild numbers conjecture implies the
weak $3x+1$ conjecture. 
As a final result we show that these conjectures
completely characterize 
the  structure of the wild semigroup,
which turns out to be quite tame.

Based on some of the
results derived here, the weak $3x+1$
conjecture and wild numbers conjecture 
were subsequently proved by David Applegate
and the author in \cite{AL04}.
In section 4  we indicate some
features of the proof 
and formulate some open problems 
about Wooley numbers, which  remain
mysterious.  

The terms  
``wild semigroup'' and ``wild number'' 
were suggested by  the novel {\it The Wild Numbers} by
Philibert Schogt \cite{Sch98}. 
The novel chronicles the efforts of a
mathematics  professor to solve 
the  (fictitious) ``Beauregard Wild Numbers Problem,''
while dealing with the ups and downs of life in
a university mathematics department.
The semigroup problem posed here has 
some  striking resemblances 
to the information given about the
Beaureguard wild numbers problem.
Beaureguard wild numbers are described in the novel as  
certain integers produced at the end of a sequence of
elementary operations  that involve 
noninteger rationals at the
intermediate steps. Here the  
semigroup products of $\sW$
giving an element of $\sW(\ZZ)$ 
generally consist of rationals
whose partial products  
typically  become integers only at
the last step. The novel also states \cite[pp.34, 37]{Sch98} that 
$2$ is a Beaureguard wild number but $3$ is not,
and  that $67$ and $4769$ are Beaureguard wild numbers.
The wild numbers defeined here reproduce nearly
all this empirical data. 
(The one exception is 
$4769 = 169\cdot 253$, which 
belongs to the wild integer semigroup $\sW(\ZZ)$
but is not a wild number as we define it.
Perhaps the novel has a misprint   
for $4759$ or $4789$ or $4967$, all primes.)
The Beaureguard
wild numbers problem is to decide whether  there
are infinitely many wild numbers (\cite[p.35]{Sch98}).
The terms  ``Wooley semigroup'' and
``Wooley numbers'' are named 
after Trevor D. Wooley, in honor of his work in
related areas of number theory (for example, \cite{BW98}).

Aside from the definitions used in this paper, 
there  have been other definitions proposed for  ``wild numbers''
that possess some of the properties indicated in the foregoing
discussion. We refer to  
sequence A58883 in the Encyclopedia of Integer Sequences
maintained by Neil Sloane \cite{Slo95}, and six
versions of ``pseudo-wild numbers'' cited there.
The Beaureguard wild numbers problem in Schogt's  novel seems
to involve iteration, which is not
directly present in our semigroup problem.  Some 
iteration problems with a similar flavor to
the wild numbers problem come from the  
``approximate multiplication'' maps
studied in Lagarias and Sloane \cite{LS04}.
A typical example is  the map
$f(x) = \frac{4}{3}\lceil x\rceil$.
The  question studied in \cite{LS04} asks
whether it is true that, for each positive integer $n$, some iterate 
$f^{(N)}(n)$ is again an integer. This  
iteration problem thus produces
a sequence of noninteger rational numbers terminating
in an integer. It is currently unsolved and seems likely
to be difficult.

%
%
%
%

\section{WOOLEY NUMBERS}
\hsp
We show that the Wooley semigroup $\sW_0$ is a
recursive semigroup. 

\begin{theorem}~\label{th21}
There is
an effectively computable procedure
that for any given positive rational $r$
determines whether
or not it belongs to the Wooley semigroup $\sW_0$,
and if it does, exhibits it as a product of generators.
\end{theorem}

\paragraph{}{\it Proof.}
We cannot represent $r$ unless it is a positive
rational number having an odd denominator (in lowest terms).
Let $g(n) = \frac{3n+2}{2n+1}$
denote the $n$th generator of the semigroup
$\sW_0$, and suppose that
$r = \prod_{i=1}^m g(n_i)$ with
$n_1 \le n_2 \le \cdots \le n_m.$
We first  bound $m$ above. 
In fact, since
$g(n) > \frac{3}{2}$ for each $n$,  we must have
$r > (\frac{3}{2})^m$, which delivers an upper bound for  $m$.

Now let $m$ be fixed. We find an  upper bound for $n_1$.
We have $r > (\frac{3}{2})^m$, so  $r=(\frac{3}{2}+\epsilon)^m$
with $\epsilon= r^{1/m} - \frac{3}{2} > 0$. We claim that
$n_1 \le 1/\epsilon.$ If not, then  
$$
g(n_1) = \frac{3n_1+2}{2n_1+1} = \frac{3}{2} + \frac{1/2}{2n_1+1}
<   \frac{3}{2} + \epsilon.
$$
Since $g(n)$ is a decreasing function of $n$, we would have
$$
r = \prod_{i=1}^m g(n_i) \le g(n_1)^m < (\frac{3}{2} + \epsilon)^m = r,
$$
a contradiction that  proves the claim. 

Once $n_1$ is chosen,  we can divide out
$g(n_1)$ to create a new problem of the same kind with
a smaller value $r' = r (g(n_1))^{-1} < \frac{2}{3}r$,
where we ask for a representation using a product of exactly
$m-1$ generators. We then show that there a finite
set of choices for $n_2$, obtaining in the process an explicit upper
bound for $n_2$ as a function of $r, m ,$ and $n_1$. 
Proceeding by
induction on $m$, we discover that the total allowable set of choices
is finite, with an effectively computable upper bound. 
Searching all of them yields either a relation certifying that
$r$ belongs to $\sW_0$ or a proof that $r$ does not belong to $\sW_0$.
$~~~\bsq$

We can carry out  this procedure in the simplest cases.

\noindent {\bf Example 2.2.} {\it
The integers $5$ and $10$ are not Wooley integers,
but  $20$ is a Wooley integer. As a consequence,  $20$ is 
a Wooley number.} \\

\noindent {\it Proof.}
Suppose that $5$ were a product of generators of $\sW_0$.
Since $2$ cannot be cancelled from the numerator of any product or
$3$ from its denominator, any representation of $5$
could not  use the generators $g(0)= \frac{2}{1}, 
g(1)= \frac{5}{3}$, or $g(2)=\frac{8}{5}$.
Any product of three of the remaining generators is
no larger than 
$g(3)^3=(\frac{11}{7})^3 = \frac{1331}{243} <5$,
 so any representation of $5$ necessarily includes at  
least four factors from
the generating set of $\sW_0$. However,
 any such product is larger than 
$(\frac{3}{2})^4 = \frac{81}{16} > 5$, a contradiction.

Suppose that $10$ were a product of generators of $\sW_0$.
Any representation of $10$ would use at most five generators, since
$(\frac{3}{2})^6 = \frac{729}{64} >10.$
A representation of $10$ could not use the generator
$\frac{2}{1}$, for if it did this fraction could be
removed, yielding a representation of $5$, a contradiction.
Also $\frac{5}{3}$ and $\frac{8}{5}$
could not arise as factors  for the same reasons, so 
the fraction of largest size that could appear in any product is again
$\frac{11}{7}.$ However $(\frac{11}{7})^5 < 10$, so there can be no
such representation.

The number $20$ can be expressed as follows:
$$
20 = g(3)^2 \cdot g(5)\cdot g(8) \cdot  g(27) \cdot g(32) \cdot  g(41) 
= (\frac{11}{7})^2 (\frac{17}{11}) ( \frac{26}{17}) (\frac{83}{55})
(\frac{98}{65})(\frac{125}{83}). 
$$
This confirms that $20$ belongs to $\sW_0$, making it a Wooley integer.
To see that $20$ is a Wooley number, note that if 
it is not irreducible, then $20=n_1 n_2$, where $n_1$ and $n_2$  
belong  to  $\sW_0(\ZZ)$. 
At least one  of $n_1$ or $n_2$ would then be divisible by $5$, 
but the only possibilities are $5$ and $10$, which have already been
ruled out.
$~~~\bsq$

The algorithm of Theorem~\ref{th21} appears
to require at least exponential time. However one can find  Wooley
integers by less exhaustive methods. Table 1 presents 
additional Wooley integers with 
identities certifying their membership in $\sW_0(\ZZ)$
for certain numbers of the form $2^k p$, where $p$ is prime
such that $5 \le p < 50$. 
These identities were found by  Allan Wilks via computer search.
Wilks's search used certain heuristics, and
did not decide whether these products give
the minimal power of $2$ possible. As a result we can say only
that the entries of the table are Wooley integers,
not necessarily Wooley numbers.

A commutative semigroup with or without unit $1$  is said to
be {\it a free commutative semigroup}
if every element of the semigroup except $1$ 
can be  factored uniquely  (up to ordering of the factors)
into a product of irreducible elements.
Many such semigroups arise in
number theory (see  Knopfmacher \cite{Kn90}).
We raise the question of whether the 
Wooley integer semigroup $\sW_0(\ZZ)$
has such unique factorization. It seems possible
that the answer will be  negative. To show this it would suffice
to find a Wooley number that contained two distinct odd prime
factors. A suggestive example is provided by 
$2^6 \cdot 31 \cdot 41$, which  is a Wooley integer 
expressible in terms of generators by
$$
2^6 \cdot 31 \cdot 41= g(423) \cdot (2^2 \cdot 7) (2^2 \cdot 11)^2,
$$
since  $g(423) = 1271/847= (31 \cdot 41)/(7 \cdot 11^2)$,
and both $2^2 \cdot 7$ and $2^2 \cdot 11$ belong to $\sW_0(\ZZ)$
according to  Table 1. 
 If there were a Wooley number of the
form $2^c\cdot 31 \cdot 41$, then the semigroup $\sW_0(\ZZ)$
would not be a  free commutative semigroup, because there would
exist four irreducible elements of $\sW_0(\ZZ)$ --
$ 2, 2^a \cdot 31, 2^b \cdot 41,
2^c \cdot 31 \cdot 41$ --leading to  
a nonunique factorization of $M= 2^{a+b+c}\cdot 31 \cdot 41$.
Such a  Wooley number will exist
unless there are Wooley integers of
form $2^a\cdot 31$ and $2^b \cdot 41$ with $a+b \le 6$, and
this possibility can be tested algorithmically, according
to Theorem~\ref{th21}.

\begin{table}
\begin{tabular}{|rrl|}
\hline
  $2^2 \cdot 5$ & $ =$ & $ 
(\frac{11}{7})^2 \cdot \frac{17}{11} \cdot \frac{26}{17} 
\cdot \frac{83}{55} \cdot \frac{98}{65} \cdot \frac{125}{83} $  \\
& $ =$ &
$g(3)^2 \cdot g(5)\cdot g(8) \cdot g(27) \cdot g(32) \cdot g(41) $  \\
& &\\

  $2^2 \cdot 7$ & $ =$ &  $\frac{11}{7}    \cdot \frac{26}{17}    
\cdot \frac{35}{23} 
\cdot \frac{215}{143} 
\cdot \frac{299}{199} \cdot  \frac{323}{215} 
\cdot \frac{371}{247} \cdot  \frac{398}{265}    $ \\
& $ =$ & 
$g(3) \cdot g(8) \cdot g(11) \cdot g(71) \cdot g(99) \cdot g(107) 
\cdot g(123) \cdot g(132)$ \\
& &\\
 
 $2^2 \cdot 11$ & $ = $ &
$(\frac{11}{7}) ^2 \cdot \frac{26}{17} \cdot \frac{35}{23}  \cdot 
\frac{215}{143} \cdot \frac{299}{199} \cdot \frac{323}{215} 
\cdot  \frac{371}{247}    \cdot \frac{398}{265} $ \\
& $ = $ &
$g(3)^2 \cdot g(8) \cdot g(11)  \cdot g(71) \cdot g(99) \cdot g(107)
\cdot g(123) \cdot g(132)$ \\
& &\\

  $2^3 \cdot 13$  & $ =$ &
$(\frac{11}{7})^2 \cdot (\frac{17}{11})^3 \cdot (\frac{26}{17})^2 
\cdot \frac{35}{23} \cdot \frac{215}{143} \cdot \frac{299}{199} \cdot 
\frac{323}{215}    \cdot  \frac{371}{247}     \cdot \frac{398}{265}  $   \\
& $ =$ &
$g(3)^2 \cdot g(5)^3 \cdot g(8)^2 \cdot g(11) \cdot g(71) \cdot g(99) 
g(107) \cdot g(123) \cdot g(132)$   \\
& &\\

  $2^2 \cdot 17$  & $ =$ &
$(\frac{11}{7})^2 \cdot \frac{17}{11} \cdot \frac{26}{17}
\cdot \frac{83}{55} \cdot \frac{98}{65} \cdot \frac{125}{83}
\cdot \frac{143}{95} \cdot \frac{215}{143} \cdot \frac{323}{215} $ \\
& $ =$ &
$g(3)^2 \cdot g(5) \cdot g(8) \cdot g(27) \cdot g(32) \cdot g(41) 
\cdot g(47) \cdot g(71) \cdot g(107)$ \\
& &\\

  $2^5 \cdot 19$  & $ =$ &
$(\frac{11}{7})^4 \cdot (\frac{17}{11})^2 \cdot (\frac{26}{17})^2
\cdot \frac{38}{25} \cdot (\frac{83}{55})^2 \cdot (\frac{98}{65})^2
\cdot (\frac{125}{83})^2 $\\
& $ =$ &
$g(3)^4 \cdot g(5)^2 \cdot g(8)^2 \cdot g(12)  \cdot g(27)^2 \cdot g(32)^2 
\cdot g(41)^2$ \\
& &\\

  $2^5 \cdot 23$  & $ =$ &
$ \frac{11}{7}\cdot \frac{26}{17} \cdot \frac{35}{23} 
\cdot \frac{47}{31} \cdot \frac{137}{91} \cdot  \frac{155}{103} \cdot 
\frac{206}{137} \cdot \frac{215}{143} \cdot (\frac{299}{199})^2
\cdot \frac{323}{215} \cdot \frac{353}{235} \cdot 
\frac{371}{247} \cdot (\frac{398}{265})^2 \cdot \frac{530}{353}$ \\
& $ =$ &
$g(3) \cdot g(8) \cdot g(11)  
\cdot g(15) \cdot g(45) \cdot  g(51) \cdot g(68) \cdot g(71)$  \\
& & ~~~~~$\cdot g(99)^2
\cdot g(107) \cdot g(117) \cdot g(123) \cdot g(132)^2 \cdot g(176)$ \\
& &\\

  $2^5 \cdot 29$  & $ =$ &
$(\frac{11}{7})^4 \cdot (\frac{17}{11})^2 \cdot (\frac{26}{17})^2 
\cdot \frac{29}{19} 
\cdot \frac{38}{25} 
\cdot  (\frac{83}{55})^2 \cdot (\frac{98}{65})^2 \cdot (\frac{125}{83})^2 $\\
& $ =$ &
$g(3)^4 \cdot g(5)^2 \cdot g(8)^2 \cdot g(9) \cdot g(12) 
\cdot g(27)^2 \cdot g(32)^2 \cdot g(41)^2 $ \\
& &\\

  $2^{11} \cdot 31$  & $ =$ &
$(\frac{11}{7})^6 \cdot (\frac{17}{11})^3 \cdot \frac{29}{19} \cdot
\frac{38}{25} \cdot \frac{62}{41} \cdot (\frac{83}{55})^3 \cdot
(\frac{98}{65})^3  \cdot (\frac{125}{83})^3 \cdot \frac{164}{109} \cdot
\frac{218}{145}$\\
& $ =$ &
$g(3)^6 \cdot g(5)^3 \cdot g(8)^3 \cdot g(9) \cdot g(12) \cdot g(20)$ \\
& & ~~~~~$\cdot g(27)^3 \cdot g(32)^3 \cdot g(41)^3 \cdot g(54) 
\cdot g(72) $ \\
& &\\

$2^5 \cdot 37$  & $ =$ &
$(\frac{11}{7})^2 \cdot (\frac{26}{17})^2  \cdot (\frac{35}{23})^2 \cdot
\frac{74}{49} \cdot (\frac{215}{143})^2 \cdot (\frac{299}{199})^2 \cdot
(\frac{323}{215})^2 \cdot (\frac{371}{247})^2 \cdot (\frac{398}{265})^2$\\
& $ =$ &
$g(3)^2 \cdot g(8)^2 \cdot g(11)^2 \cdot g(24)
\cdot g(71)^2 \cdot g(99)^2 \cdot g(107)^2  \cdot g(123)^2 \cdot  g(132)^2$\\
& &\\

  $2^{10} \cdot 41$ & $ =$ &
$(\frac{11}{7})^6 \cdot (\frac{17}{11})^3 \cdot (\frac{26}{17})^3  \cdot
\frac{29}{19} \cdot \frac{38}{25} \cdot (\frac{83}{55})^3 \cdot
(\frac{98}{65})^3  \cdot (\frac{125}{83})^3 \cdot \frac{164}{109} \cdot
\frac{218}{145}$\\
& $ =$ &
$g(3)^6 \cdot g(5)^3 \cdot g(8)^3 \cdot g(9) \cdot g(12) 
\cdot g(27)^3 \cdot g(32)^3 \cdot g(41)^3 \cdot g(54) \cdot g(72) $ \\
& &\\

$2^{11} \cdot 43$ & $ =$ &
$ (\frac{11}{7})^5 \cdot (\frac{17}{11})^2 \cdot(\frac{26}{17})^3 
\cdot \frac{29}{19} \cdot \frac{35}{23} \cdot \frac{38}{25} \cdot
(\frac{83}{55})^2 \cdot (\frac{98}{65})^2 \cdot (\frac{125}{87})^2 
\cdot \frac{215}{143}$  \\
& &~~~~~  $ \cdot \frac{299}{199}    \cdot \frac{305}{203} \cdot 
\frac{323}{215} \cdot \frac{344}{229} \cdot \frac{371}{247} 
\cdot \frac{398}{265} \cdot \frac{458}{305}$ \\
& $ =$ &
$ g(3)^5 \cdot g(5)^2 \cdot g(8)^3 \cdot g(9) \cdot g(11) \cdot g(12) \cdot
g(27)^2 \cdot g(32)^2 \cdot g(41)^2 \cdot g(71)$  \\
& &~~~~~  $ \cdot g(99) \cdot g(101) \cdot g(107) \cdot g(114) \cdot g(123) 
\cdot g(132) \cdot g(152)$ \\
& &\\
  $2^{11} \cdot 47$  & $ =$ &
$(\frac{11}{7})^6 \cdot (\frac{17}{11})^3 \cdot \frac{29}{19} \cdot
\frac{38}{25} \cdot \frac{47}{31} \cdot \frac{62}{41} \cdot 
(\frac{83}{55})^3 \cdot
(\frac{98}{65})^3  \cdot (\frac{125}{83})^3 \cdot \frac{164}{109} \cdot
\frac{218}{145}$\\
& $ =$ &
$g(3)^6 \cdot g(5)^3 \cdot g(8)^3 \cdot g(9) \cdot g(12) \cdot g(15) 
\cdot g(20)$ \\
& & ~~~~~$\cdot g(27)^3 \cdot g(32)^3 \cdot g(41)^3 \cdot g(54) 
\cdot g(72) $ \\
& &\\

\hline
\end{tabular}
\caption{Members of the Wooley integer semigroup $\sW_0(\ZZ)$.}
\label{table1}
\end{table}

%
%
%
\section{WILD NUMBERS} 
\hsp
We begin by showing that there are infinitely
many wild numbers.

\begin{theorem}~\label{th31}
The semigroup of wild integers
contains infinitely many irreducible elements (i.e., there are
infinitely many wild numbers.)
\end{theorem}

\paragraph{}{\it Proof.}
For $n = \frac{5^k -1} {2}$ we have
$g(n) = \frac{\frac{1}{2}(3 \cdot 5^k + 1)}{5^k}.$
Example 2.2 shows  that $2^2 \cdot 5$ is a Wooley  number, 
which implies that $5$ is a wild integer, and it is a
wild number since it is prime. We conclude that
$$
h(k) :=\frac{1}{2}(3 \cdot 5^k + 1) = g(n) \cdot 5^k
$$
belongs to  $\sW(\ZZ)$ for each positive integer $k$.

The sequence $\{h(k): k \ge 1\}$ satisfies a
homogeneous second-order linear recurrence, namely,
$$
h(k) = 6 h(k-1) - 5 h(k-2).
$$
This sequence 
 is nondegenerate in the sense of Ward \cite{Wa54}
(i.e., it does not satisfy a first-order linear recurrence).
Accordingly, by  the main result of Ward \cite{Wa54} the sequence $\{h(k)\}$
contains an infinite number of distinct prime divisors
(i.e., the set $D$ of primes $p$ that divide $h(k)$ for
at least one $k$ is infinite). 

We now argue by contradiction that 
 $\sW(\ZZ)$ contains infinitely many irreducible elements.
If not, there would exist  some prime $p$
in the infinite set $D$ that did not divide 
any irreducible element. This prime $p$ divides
some $h(k)$, which belongs to $\sW(\ZZ)$,
so there exists a smallest 
element $m$ of $\sW(\ZZ)$ that is divisible by $p$.
This element $m$ is necessarily irreducible,
for if not there would be a smaller integer in $\sW(\ZZ)$
divisible by $p$. This gives a contradiction.
$~~~\bsq$

We can obtain  much stronger results about
the structure of the wild integer semigroup if we assume
the truth of the weak $3x+1$ conjecture. 

\begin{theorem}~\label{th32}
Suppose that  the weak $3x+1$ conjecture holds.
Then the wild integer  semigroup  $\sW(\ZZ)$ is 
a free commutative semigroup whose set of generators $\sP$ consists
entirely of primes. In other words, all wild numbers are
prime numbers.
\end{theorem}

\paragraph{}{\it Proof.} 
The  semigroup $\sW(\ZZ)$ contains all powers of $2$.
Also, since  $2$ is invertible in $\sW$, 
it can be cancelled from all other generators, which
therefore must be
odd integers. However, the weak $3x+1$ conjecture says
that if $n$ belongs to $\sW$, then so does $\frac{n}{k}$ for any
positive integer $k$. Thus, if a composite number $n$ 
lies in $\sW(\ZZ)$,
so do all of its prime divisors. It follows that all
generators of $\sW(\ZZ)$ are primes. The semigroup $\sW(\ZZ)$
is now a free commutative  semigroup as a consequence of unique prime
factorization of integers. 
$~~~\bsq$

In general, we can certify that a given prime number $p$
is a wild number by finding  some $j$ such that $2^j p$  is a
Wooley number. For example,  $67$ is a wild number since
$2^{12} \cdot 67$ is a Wooley number. The latter assertion
is a consequence of the identity
$$
\frac{2^5 \cdot67}{5 \cdot 37} = g(29) \cdot g(44) \cdot g(69) \cdot
 g(78) \cdot g(92) \cdot g(104)
$$
and the 
fact,  established
earlier, that $2^2 \cdot 5$ and $2^5 \cdot 37$ are Wooley  numbers.

We next show that the weak $3x+1$ conjecture 
implies the 
wild numbers conjecture.

\begin{theorem}~\label{th33}
If  the weak $3x+1$ conjecture holds, then 
the wild numbers conjecture is true.
\end{theorem}

\paragraph{Proof.}
We prove the wild numbers conjecture by induction on
the $n$th prime, call it $q$. The 
the induction hypothesis asserts that all smaller primes
except $3$ belong to the wild integer semigroup.
We call an integer
{\it $Y$-smooth} if all its prime factors are strictly
smaller than $Y$. Because the wild integers form
a semigroup, the induction hypothesis tells us
that all $q$-smooth numbers not divisible by $3$ 
are wild integers.
In particular, all  integers smaller than $q$ 
and not divisible by $3$ are
wild integers. 
To complete the induction step it suffices to  show that
$q$ is a wild integer. If so, its primality
guarantees that it is irreducible, 
hence that it is  a wild number.

To show $q$ is a wild integer it  suffices to find
some multiple $mq$ that is a wild integer, for
the weak $3x+1$ conjecture implies that
 $1/m$ belongs to the wild semigroup $\sW$,
making  $q =  (mq)/m$  a wild integer.
We wish to find $mq$
of the form
$mq= 3n+2$ such that  $2n+1$ is a $q$-smooth number
not divisible by $3$.
If so, then $mq = \frac{3n+2}{2n+1} (2n+1)$ will
belong to  $\sW$,
and the desired result will follow. The congruence
restriction  puts $m$ in a certain residue class modulo $3$,
and by imposing a condition modulo $9$ we can guarantee
that $2n+1 \not\equiv 0 ~(\bmod~3)$. The resulting  
integers $2n+1$ fall into  an arithmetic progression
of numbers congruent to $r$ modulo $6q$, with $gcd(r, 6q)=1$, 
and we arrive at a special case of the well-studied
arithmetic question of finding ``smooth numbers'' in
an arithmetic progression.

We recall general facts on the distribution
of ``smooth numbers'' up to X, namely, those numbers below 
$X$ having all prime factors smaller than a given  bound $Y$
(see,  for example,  Hildebrand and Tenenbaum \cite{HT93}).
Smooth numbers with an appropriate choice of $Y$ 
play an important role in the 
design and performance of the   fastest known
algorithms for factoring large numbers
(see  Pomerance \cite{Po95}).
It is known that the number of integers  smaller than
$X$ that have all their prime factors below a cutoff value 
$Y=X^{\alpha}$ for any fixed $\alpha$
have asymptotically a positive density
$\rho(\alpha)Y$, where $\rho(u)$ is a strictly positive function, 
the {\em Dickman function}. This function is  given by the
solution to a certain difference-differential equation 
and   $\rho(u) \approx u^{-u}$.
(Here $u = (\log Y)/(\log X)$.)
Results of Balog and Pomerance \cite{BP92} 
carry these bounds over 
to count the number of
$Y$-smooth numbers in arithmetic progressions modulo N, and their results
give an asymptotic formula valid for $Y$ over
a large range. In particular,
choosing $X \approx q^2$, $N = 6q$, and $Y = q$,
which corresponds to  $\alpha= 1/2$, 
one can deduce from their results 
 that the number of $q$-smooth integers in the first $q$ terms of the
arithmetic progression $r$ modulo $6q$ is nonzero whenever $q > C_0$,
for some constant $C_0$ that is,  
in principle, computable. However, $C_0$ is not easy to
compute, nor is it likely to be small.
 
To complete our argument we need only demonstrate
the existence of a single $q$-smooth number in the given arithmetic
progression   $r$ modulo $6q$. This
permits us to sidestep the results of Balog and Pomerance 
and to use instead a direct combinatorial argument. 
It rests on the observation  
that if more than
half the  invertible residue classes  modulo $N$
contain $Y$-smooth numbers, then 
(by the pigeonhole principle) every invertible
residue class $r$  modulo $N$
occurs as a  product of two
of these residue classes, and consequently contains a $Y$-smooth number
that is the product of $Y$-smooth numbers from these classes.
In more detail, let $\Sigma$
denote the set of invertible residue classes  $s$ modulo $N$
that contain a $Y$-smooth integer.
Suppose that  $r$  is an arbitrary invertible residue class modulo $N$.
We now define $\Sigma'$ to consist of 
those residue classes $s'$ modulo $N$ given by  
$$
s' \equiv r \cdot s^{-1} (\bmod~N),
$$
where $s$ belongs to $\Sigma$. Certainly $|\Sigma'|= |\Sigma|$
and since  $|\Sigma|$ and $|\Sigma'|$
exceed half the invertible residue classes,
there must be some $s'$ in $\Sigma \cap \Sigma'$.
Now $r \equiv s s'~(\bmod~N)$, and 
taking $S$ and $S'$ to be $Y$-smooth numbers in the classes
$s$ and $s'$ modulo $N$, respectively, we find that $S\cdot S'$
is a $Y$-smooth number in the class $r$ modulo $N$.

In our case we have $N=6q$, which has  $\phi(N)=2(q-1)$
invertible residue classes, and all of these residue classes
consist of numbers not divisible by $3$. It suffices to
demonstrate that  more than $q-1$ of these invertible
classes contain $q$-smooth numbers.
We show that the set of such residue
classes whose least positive residue is smooth exceeds $q-1$
for all sufficiently large $q$.
Now every 
integer less than  $6q$ and relatively prime
to $6q$ has all its prime factors smaller than $q$, 
except for primes $p'$ with $q < p' < 6q$
and integers of the form $5p'$ with $q \le p' < \frac{6}{5}q$.
Since the number of primes below
$x$ is $O (\frac{x}{\log x})$,
there are at most $O(\frac{q}{\log q})$ such integers, hence
at least $2(q-1) - O (\frac{q}{\log q})$ invertible classes
have their least residue being $q$-smooth. For large $q$ this 
gives the result, and by obtaining explicit numerical
bounds for the remainder term it is possible to prove that 
for $q > 10^4$ more than  half the invertible residue classes 
modulo $6q$ are $q$-smooth.
(Such bounds are explicitly derived in \cite{AL04}.)

Now we can prove the wild numbers conjecture by
induction, under the assumption that the weak
$3x+1$ conjecture is valid, with the base case
consisting of checking all $q$ such that $q < 10^4.$ 
The base case can be checked by computer. 
In fact, it suffices to use Table 1 when 
$q < 50$ and when $50 < q < 10^4$, to 
have the computer directly find  a smooth number in a suitable arithmetic
progression modulo $6q$.
$~~~\bsq$

Since the $3x+1$ conjecture appears to be true, Theorem~\ref{th33}
provides a powerful argument in favor of the
wild numbers conjecture. On the other hand, we
have a converse implication: 

\begin{theorem}~\label{th34}
If the wild numbers conjecture is true, 
then the weak $3x+1$ conjecture holds.
\end{theorem}

\noindent {\it Proof.}
This implication is proved in a fashion similar to the argument of
Theorem~\ref{th33}. 
We proceed by induction on the $n$th prime $q$, assuming
that all primes below $q$ (including $3$)
belong to the inverse semigroup $\sW^{-1}$.
We now consider multiples $mq$, where 
$m \equiv 1~(\bmod~ 6)$. The wild numbers conjecture
implies that all such integers are in the wild
semigroup $\sW$. Writing $mq= 2n+1$, 
we look for a case in which  $3n+2$ is a $q$-smooth number.
Expressing $m$ as 
$m=6k+1$, we have 
$$
3n+2 = 9kq + \frac{3q+1}{2},
$$
 which is an
arithmetic progression  modulo $9q$. 
As in the earlier result,
 it suffices to show that
for all sufficiently large $q$ more than half of the 
invertible residue classes modulo $9q$ in
the interval $[1, q-1]$ have least positive residues that are 
$q$-smooth numbers,
which then implies that
each arithmetic progression  for an invertible
residue class modulo $9q$ 
contains a $q$-smooth integer smaller than $81q^2$.
This holds when $q > 10^5.$ Since the $3x+1$
conjecture has been verified up to $10^5$, the
base case of the induction is already done.
$~~~\bsq$.

Our final result points out that the truth of the 
weak $3x+1$ conjecture completely determines
the structure of the wild semigroup $\sW$.

\begin{theorem}~\label{th35}
If the weak $3x+1$ conjecture 
is true, then the wild semigroup $\sW$
consists of all positive rational numbers  $a/b$ 
with $gcd(a, 3b)= 1$.
\end{theorem}

\noindent{\it Proof.} 
The weak $3x+1$  conjecture implies that $\sW$ contains
all fractions $1/p$, where $p$ is prime.
By Theorem~\ref{th33} this conjecture
 implies the wild numbers conjecture, which
ensures  that $\sW$  contains
all primes $p$ different from $3$. We observed earlier that
any rational member $r= a/b$ of $ \sW$ in lowest
terms has numerator $a$ relatively prime to $3$. This gives
the result.
$~~~\bsq$

Theorem ~\ref{th35}
provides  a simple  effective decision procedure for membership of
a given rational number $r$ in the wild semigroup  $\sW$, 
provided that the weak $3x+1$ conjecture is proved.
%
%
%
\section{CONCLUDING REMARKS} 
\hsp
The results of section 3 demonstrate that the wild numbers conjecture
and the weak $3x+1$ conjecture are intertwined: each
is implied by the other. 
David Applegate and the author \cite{AL04} have
recently been able 
to prove both conjectures simultaneously, via a bootstrap
induction procedure that uses the truth of one of the conjectures
on an interval to extend the truth of the other to a
larger interval, and vice versa. The 
argument of Theorem~\ref{th33}
(respectively,  Theorem ~\ref{th35}) provides a way to
extend the truth of the conjecture in  one direction,
provided it holds on a sufficient initial interval  in the other. 
We  do not, however, know a 
way to use the arguments of these theorems
 simultaneously in  both directions.
Fortunately, there 
is another  systematic way to find representations of 
many integers $n$ in the inverse semigroup $\sW^{-1}$,
which is  to iterate the $3x+1$ map starting with $n$.
The argument in  \cite{AL04}  takes advantage of this fact
in constructing
the ``other'' direction of the bootstrap induction.
There is an apparent asymmetry in the two directions, in
that we are not aware of any dynamical system associated
with the wild semigroup that  produces relations generating  
the integers in the wild integer semigroup $\sW(\ZZ)$
that is analogous to the use of the $3x+1$ iteration in the
inverse semigroup $\sW^{-1}$

The original motivation for studying the wild semigroup
came from the weak $3x+1$ conjecture, but the
Wooley semigroup that arose in the process seems
interesting in its own right. 
The Wooley integer semigroup $\sW_0(\ZZ)$ 
appears to be a  more complicated object than the wild
integer semigroup $\sW(\ZZ)$. 
There remain many open questions about Wooley integers.
One question already raised 
in Section 2 asks whether Wooley integers have
unique factorization into irreducibles, i.e., whether 
the Wooley integer semigroup is
a free commutative semigroup. 
A second  question  concerns, for each
prime $p$, the behavior of the minimal 
power $e(p)$ necessary to place  $2^{e(p)} p$ in
the Wooley integer semigroup. It seems plausible that $e(p)$
is unbounded. The truth of the  wild numbers conjecture
implies that  each number
$e(p)$ is finite, so in view of its proof
in \cite{AL04}, this question is well posed.
A third question asks:
How does the counting function of the
Wooley  integers below $x$ grow as $x \to \infty$? 

The  wild numbers conjecture 
was named after the (fictitious) mathematical 
problem in Philibert Schogt's novel  {\it The Wild Numbers.}
In the novel the  Beauregard Wild Numbers Problem 
 was presented as a famous unsolved problem,
with a long and illustrious history.  
Its real-life namesake  fails to have either of these attributes.
Indeed it  has a short history,
and  the problem of the infinitude of wild numbers 
was settled by Theorem~\ref{th31}. 
Nevertheless,  our terminology seems fitting, for
 the novel asserts there is ``a fundamental
relationship between wild numbers and prime numbers''
\cite[p. 36]{Sch98}, and 
the  wild numbers of this paper coincide with 
the prime numbers, excluding $3$.
Understanding the 
behavior of prime numbers is 
one of the great mysteries of mathematics,
with a history as long and 
impressive as one could hope for;
see Derbyshire \cite{De03} or du Sautoy \cite{DS03}
for recent accounts.

\paragraph{ACKNOWLEDGMENTS.} Most of this
work was done while I was at A.T. \& T. Labs-Research,
whom I thank for support.  I am indebted to Allan Wilks
for computing the data in Table 1, and to Jim Reeds for
 recommending  the novel {\it The Wild Numbers}.
  Finally I thank the reviewers for helpful comments.

%
%
%


\end{document}